\documentclass{article}
\usepackage{amsmath}
\usepackage{amssymb}
\usepackage{amscd}

\begin{document}

\title{ON THE STRUCTURE OF BOL ALGEBRAS}
\author{Thomas B. Bouetou\\ \'Ecole Nationale Sup\'erieure Polytechnique,\\
B.P.~8390 Yaound\'e, Cameroun\\
e-mail:tbouetou@polytech.uninet.cm}

\maketitle \begin{quote}{\bf Abstract.}
\em The fundamental ideas of the definition of solvable and semisimple 
Bol algebras are given and some related theorems.
\end{quote}

\maketitle \begin{quote}{\bf Keywords:}
\em Bol Algebras, solvable and semisimple Bol Algebras, Lie algebras, Killing-Ricci Form for Bol Algebras
\end{quote}
\maketitle \begin{quote}{\bf AMS subject classification 2000:}
 05B07, 17D05, 17D10, 17D99, 20N05, 22E60, 22E67, 
\end{quote}

\section{Introduction}

It is well known that  structure theory is a basic tool for the manipulation
of any algebraico-geometrical object, hence this study. In this paper we 
intend to introduce the notion of solvability for any 
ideal
of Bol algebra and semisimplicity. Our approach uses the notion of solvability 
introduce
in \cite{kik3} and the one of \cite{lis}. By using properties of the defined 
expression, we will be able to define the notion of the radical of a finite 
dimensional Bol algebra. Using the notion of enveloping pair on which
Bol algebras are rooted we establish the fact that if any Bol algebra is
solvable then its enveloping Lie algebra is solvable; like wise if it is simple. As
in \cite{bou}, it is shown that the Levi-Malcev theorem is not valid for
Bol algebra in the classical sense but for a class of Bol algebras called
homogeneous. Thus, Bol algebra allows deformation, but the deformation
is quite different from the deformation associated with Lie algebras.  Bol
 algebras are also significant in Mathematical Physics.

\section{Basic definition}

{\bf Definition 1.} Any vector space $V$ over the field of characteristic $0$  with the operations 
$$\xi,\eta \to \xi\cdot\eta\in V,$$
$$\xi,\eta,\zeta\to (\xi,\eta,\zeta)\in V \quad (\xi,\eta,\zeta\in
V)$$ and identities $$ \xi\cdot\xi=0,\quad
(\xi,\eta,\zeta)+(\eta,\zeta,\xi)+(\zeta,\xi,\eta)=0,$$

$$(\xi,\eta,\zeta)
\chi-(\xi,\eta,\chi)\zeta+(\zeta,\chi,\xi\cdot\eta)-
(\xi,\eta,\zeta\cdot\chi)+\xi\eta\cdot\zeta\chi=0,$$
$$(\xi,\eta,(\zeta,\chi,\omega))=((\xi,\eta,\zeta),\chi,\omega)+(\zeta
,(\xi,\eta,\chi),\omega)+(\zeta,\chi,(\xi,\eta,\omega))$$
is called a Bol algebra.

We will note that any Bol algebra can be realized as the tangent 
algebra to a  Bol Loop with the left Bol identity, and they allow embedding in Lie
algebras. This can be expressed this way;
let $ (G, \Delta, e) $- be a local Lie group,  $ H $-  one of its subgroups, 
and let us  denote the  corresponding Lie algebra and subalgebra by
$ \mathfrak{G}$ and $ \mathfrak{h}$. Consider a vector subspace 
$\mathfrak{B}$ such that 

$ \mathfrak{G}=\mathfrak{h}\dotplus \mathfrak{B}$.
 Let $ \Pi:G\longrightarrow G\setminus H $ be the canonical projection and let
$ \Psi$ be the restriction of mappings composition $ \Pi \circ \exp$, 
to $ \mathfrak{B}$. Then there exists such a neighborhood $ \mathcal{U}$ of 
the
point $ O $ in $\mathfrak{B}$ such that $\Psi$ maps it diffeomorphically into the 
neighborhood $\Psi (u)$ of the coset $\Pi (e)$ in $ G \setminus H $ .\\

$$
\begin{CD}
\mathfrak{G} @<{i}<< \mathfrak{B}\\
@V{\exp}VV  @VV{\Psi}V \\
G @>{\pi}>> G\setminus H
\end{CD}
$$

by introducing a local composition law:

$$
 a \star b= \Pi_{B}(a\Delta b) 
$$ 
 
On points of local cross-section $ B=\exp \mathcal{U} $ of left coset of $G$
 mod $H$ where $ \prod_{B}=\exp \circ \Psi^{-1} \circ\Pi: G \longrightarrow B $
is the local projection on $ B $ parallel to $ H $ which puts  every element 
$ a \in B $ in correspondence so that $ g=a \Delta p $, where $  p \in H $.
We will emphasize that if any two local analytic Loops are isomorphic then 
their
corresponding Bol algebras are isomorphic.

\section{Solvable Bol algebras}

Let $\mathfrak{B}$ be a Bol algebra over a field $K$ of characteristic $0$.

{\bf Definition 2.} We say that the subspace $\mathfrak{V}$ is a subsystem of
$\mathfrak{B}$ if $(\mathfrak{V},\mathfrak{V},\mathfrak{V}) \subset\mathfrak{V}$
and an ideal if $(\mathfrak{V},\mathfrak{B},\mathfrak{B}) \subset\mathfrak{V}$ 
and $\mathfrak{V}\cdot \mathfrak{B} \subset\mathfrak{V}$.

Our idea is by following the general notion and keeping out the condition
of Bol algebras to give the definition of solvability of Bol algebras.
Following \cite{lis} we define the series:
$$ \mathfrak{V}^{(n+1)}=(\mathfrak{V}^{(n)},\mathfrak{V}^{(n)},\mathfrak{B}) \; \forall n \geq 0 $$

where $\mathfrak{V}$ is an ideal  and one can verify that:
$$ \mathfrak{V}^{(n+1)}\subseteq \mathfrak{V}^{(n)}\subseteq \mathfrak{V}^{(n-1)}\subseteq ......\subseteq\mathfrak{V}^{(0)}=\mathfrak{V}$$

  We say that the ideal $\mathfrak{V}$ is solvable if there
exist $k \geq 0$ such that  $\mathfrak{V}^{(k)}=0$. 

\; \; The given definition is developed in \cite{lis} for a trilinear operation.
 For Bol algebras we have 
a bilinear and a trilinear operation. Therefore, we will now define 
the ideal of Bol algebras, by  considering the two operations.

\; \; Let $\mathfrak{W}$ be an Ideal related to the l.t.s of Bol 
algebras $\mathfrak{B}$ and let us define the subsystem $\mathfrak{M}$ by
$\mathfrak{M}=\mathfrak{W}\cdot \mathfrak{W}+(\mathfrak{W},\mathfrak{W},\mathfrak{B})$

it is clear that $\mathfrak{M}$ is an ideal of  $\mathfrak{W}$ even normal.
We  define the relation:

$$ \mathfrak{W}= \mathfrak{W}^{(0)}$$

$$\mathfrak{W}^{(n+1)}=\mathfrak{W}^{(n)}\cdot \mathfrak{W}^{(n)}+(\mathfrak{W}^{(n)},\mathfrak{W}^{(n)},\mathfrak{B})$$

{\bf Theorem 1.} $\forall n>0$  $\mathfrak{W}^{(n)}$ is an ideal in 
 $\mathfrak{W}^{(n-1)}$.

  Proof:

for proof  of this theorem, we will use the induction on $n$.

For $n=1$ we must verify that 
$\mathfrak{W}\cdot \mathfrak{W}' \subset \mathfrak{W}'$ and
$(\mathfrak{W}', \mathfrak{W}, \mathfrak{W}) \subset \mathfrak{W}'$
where
$\mathfrak{W}'=\mathfrak{W} \cdot \mathfrak{W}+(\mathfrak{W}, \mathfrak{W}, \mathfrak{B})$. 

We have,
$\mathfrak{W}\cdot \mathfrak{W}'= \mathfrak{W}\cdot\mathfrak{W}\cdot \mathfrak{W}+\mathfrak{W}( \mathfrak{W}, \mathfrak{W}, \mathfrak{B})$
considering the identities of the definition of Bol algebras we have
$\mathfrak{W}\cdot \mathfrak{W}' \subset \mathfrak{W}'$ and

$(\mathfrak{W}', \mathfrak{W}, \mathfrak{W})=(\mathfrak{W}, \mathfrak{W},( \mathfrak{B}, \mathfrak{W}, \mathfrak{W})) \subset \mathfrak{W}'$.

We assume that  $\mathfrak{W}^{(n)} \subseteq \mathfrak{W}^{(n-1)}$ and
show that  $\mathfrak{W}^{(n+1)} \subseteq \mathfrak{W}^{(n)}$ or
 
$$\mathfrak{W}^{(n+1)}\subseteq \mathfrak{W}^{(n-1)}\cdot \mathfrak{W}^{(n-1)}+(\mathfrak{W}^{(n-1)},\mathfrak{W}^{(n-1)},\mathfrak{B})\equiv \mathfrak{W}^{(n)}.$$
We have 
$$\mathfrak{W}^{(n-1)}\cdot \mathfrak{W}^{(n)}\subseteq \mathfrak{W}^{(n-1)}\cdot \mathfrak{W}^{(n-1)} \subseteq \mathfrak{W}^{(n)}.$$
Also since  $\mathfrak{W}$ is a subset we have:
$$(\mathfrak{W}^{(n-1)},\mathfrak{W}^{(n-1)},\mathfrak{W}^{(n)})\subseteq (\mathfrak{W}^{(n-1)},\mathfrak{W}^{(n-1)},\mathfrak{B})\subseteq \mathfrak{W}^{(n)}$$
hence the result.

As in the case of Lie algebras, we can then write the following ascending
 series:

$$\mathfrak{W}^{(n+1)}\subset \mathfrak{W}^{(n)}\subset \mathfrak{W}^{(n-1)}\subset ...............\subset \mathfrak{W}^{(1)}\subset \mathfrak{W}^{(0)}=\mathfrak{W}.$$

After  conducting the examination one can state the following definition:

{\bf Definition 3.} The subsystem $\mathfrak{W}$ of a Bol algebra $\mathfrak{B}$
is an ideal of  $\mathfrak{B}$ if

$\mathfrak{W}\cdot \mathfrak{W}+(\mathfrak{W},\mathfrak{W},\mathfrak{B})\subset \mathfrak{W}$.

 The ideal $\mathfrak{W}$ of a Bol algebra $\mathfrak{B}$ is
said to be solvable if there exist  $k \geq 0$ such that $\mathfrak{W}^{(k)}=0$. 
 
We will discuss the notion of maximal ideal:

Lemma If $\mathfrak{V}$ and  $\mathfrak{W}$ are solvable ideals then so is
 $\mathfrak{V+W}$.

Proof 

 Assume that  $n$, $k$  exist, such that $V^{(n)}=0$ and  $W^{(k)}=0.$
 We also  suppose that $n=k$
$$V^{(n)}=V^{(n-1)}\cdot V^{(n-1)}+(V^{(n-1)},V^{(n-1)},\mathfrak{B})$$
$$W^{(n)}=W^{(n-1)}\cdot W^{(n-1)}+(V^{(n-1)},W^{(n-1)},\mathfrak{B})$$

We  show that  $(V+W)^{(n)}=0$ by using induction on $n$. For $n=0$ 
it's obvious. For $n=1$ we have

$(V+W)^{(1)}=(V+W)\cdot (V+W)+(V+W,V+W,\mathfrak{B})$\\
$\; \; \; \; \; \; \;\;=V\cdot V+W\cdot W+V\cdot W+W\cdot V+(V,V,\mathfrak{B})+(V,W,\mathfrak{B})+(W,V,\mathfrak{B})+(W,W,\mathfrak{B})$\\
$\; \; \; \; \; \; \; \;=V^{(1)}+W^{(1)}+V\cdot W+W\cdot V+(V,W,\mathfrak{B})+(W,V,\mathfrak{B})$\\
$\; \; \; \; \; \; \; \; \;\subseteq V^{(1)}+W^{(1)}+V\cap W.$

We assume that: 

$(V+W)^{(k)}\subseteq V^{(k)}+W^{(k)}+V\cap W$

and  prove that:

$(V+W)^{(k+1)}\subseteq V^{(k+1)}+W^{(k+1)}+V\cap W.$

We have

$(V+W)^{(k+1)}=(V+W)^{(k)}\cdot (V+W)^{(k)}+((V+W)^{(k)},(V+W)^{(k)},\mathfrak{B})$\\
$\; \; \; \; \; \; \; \; \subseteq ( V^{(k)}+W^{(k)}+V\cap W)(V+W)^{(k+1)}+( V^{(k)}+W^{(k)}+V\cap W,(V+W)^{(k)},\mathfrak{B}).$

After opening and rearranging we obtain the result.

This leads us to define the radical of Bol algebras.

{\bf Definition 4.} The radical $\mathfrak{R}$ of a Bol algebra $\mathfrak{B}$ is
the maximal and unique solvable ideal.
 A Bol algebra is said to be  semisimple if its radical is
 zero.

As they are constructed, Bol algebra can be  enveloped in Lie algebra. 
 We now discuss some facts about it using the notion of pseudo-derivation.
Following \cite{sab1} we know Bol algebras are defined from two operations
$(\cdot)$ and $( ; ,)$ verifying a series of conditions. In what follows, 
we will discuss some results from \cite{sab1}.

{\bf Definition 5.} The linear endomorphism $\prod$ of the binary-ternary
algebra $\mathfrak{B}$ with the composition law $X\cdot Y$ and $(X ;Y ,W)$ 
will 
be called  pseudo-derivation with the component $Z$ if
$$\prod (X\cdot Y)=(\prod X)\cdot Y+X\cdot (\prod Y)+(Z;X,Y)+(X\cdot Y)\cdot Z$$
$$\prod (X; Y,W)=(\prod X;Y,W)+ (X;\prod Y,W)+(X;Y,\prod Z) $$
$$ \forall X,Y,Z,W \in \mathfrak{B}$$

The  pseudo-derivation of Bol algebra $\mathfrak{B}$ form a Lie algebra with
respect to the operation of sum and multiplication by a scalar and the Lie
commutator $[\prod,\widetilde{\prod}]=\prod \widetilde{\prod}- \widetilde{\prod}\prod $.

For Bol algebra, one can verify that if  $\prod $, and $\widetilde{\prod}$ has
$Z$ and $\widetilde{Z}$ as a component respectively then  $\prod+\widetilde{\prod}$
and $\lambda \prod $ have $Z+\widetilde{Z}$ and $\lambda Z$ correspondently 
as a component but
 $[\prod,\widetilde{\prod}]$ has for a component $Z\cdot \widetilde{Z}+\prod \widetilde{Z}-\widetilde{\prod}Z$. We  denote the algebra formed by the  
pseudo-derivation by  $pder\mathfrak{B}$.
 Using this notation we  re-write the definition of Bol algebra as follows:

{\bf Definition 6.} The binary- ternary $\mathbb{R}$-algebra $\mathfrak{B}$
with the bilinear $(\cdot)$ and the trilinear $(;,)$ operation is called
a Bol algebra if:
\begin{itemize}
\item $X\cdot X=0$
\item $(X;Y,Y)=0$
\item $(X;Y,Z)+(Z;X,Y)+(Y;Z,Y)=0$
\item the endomorphism $\mathbb{D}_{X,Y} : Z\longrightarrow (Z;X,Y). $ its a pseudo-derivation with the component  $X\cdot Y$.
\end{itemize}

From \cite{sab1} it is shown that the set
 
$Pder \mathfrak{B}=\left\{(\prod,a);\prod \subset \prod p\;der\mathfrak{B}, a\; composent \;of\; \prod \right\}$

is a Lie algebra under a proper operation.

 This notion of  pseudo-derivation helps 
to 
define the enveloping Lie algebra for a Bol algebra. In \cite{sab1} it is 
proved that for any Bol algebra $\mathfrak{B}$ with the operation
$X\cdot Y$ and $(Z;X,Y)$ there is a pair of algebras
 $(\mathfrak{G},\mathfrak{h})$ such that 
$\mathfrak{G}=\mathfrak{B}\dotplus \mathfrak{h}$, and
$[\mathfrak{B},[\mathfrak{B},\mathfrak{B}]\subset \mathfrak{B}$,
$[\mathfrak{B},\mathfrak{B}]\cap \mathfrak{B}=\{0\}$ with 
$pro.\mathfrak{B}_{[X,Y]}=X\cdot Y$ and $(Z;X,Y)=[Z,[X,Y]]$ where
$X,Y,Z \in \mathfrak{B}.$  This implies that
$[X,Y]=X\cdot Y+[-X\cdot Y+(\mathfrak{D}_{X,Y},X\cdot Y)]$.

If $V$ is an ideal of Bol algebra $\mathfrak{B}$ and $\mathfrak{G}$ its  
universal enveloping Lie algebra, then any ideal $W$ generated by $V$ 
coincide with $V\dotplus [V,\mathfrak{B}]$.

We have the following theorem:

{\bf Theorem 2.} Let $V$ be a solvable ideal in $\mathfrak{B}$ then
$W=V\dotplus [V,\mathfrak{B}]$ is also a  solvable ideal in $\mathfrak{B}$.

Before proving this theorem, we state the following lemma.

{\bf Lemma 1.} If $V$ is a solvable ideal in $\mathfrak{B}$ and $\widetilde{V}$ an
extension of $V$ such that $\widetilde{V}=V+[V,\mathfrak{B}]$, then
$\widetilde{V}$  is solvable.

Proof.

We want to show that there $\exists k$ such that $ (\widetilde{V})^{(k)}=0$.

We have 
$$(V+[V,\mathfrak{B}])^{(k)}=(V+[V,\mathfrak{B}])^{(k-1)}\cdot (V+[V,\mathfrak{B}])^{(k-1)}+((V+[V,\mathfrak{B}])^{(k)},(V+[V,\mathfrak{B}])^{(k)},\mathfrak{B})$$

Here we see that the power at the left side is a function of $k$. We will then
write $p=p(k)$ therefore, we have the following lemma:

{\bf Lemma 2.} the following inclusion  holds:
$\widetilde{V}^{(p)}\subset V^{(k)}+[V^{(k)},\mathfrak{B}]$ where
$\widetilde{V}^{(0)}=\widetilde{V}$ and 
$\widetilde{V}^{(n+1)}=[\widetilde{V}^{(n)},\widetilde{V}^{(n)}]$.

Proof.

The proof of this lemma is by induction. Here we will assume that, if $k=0$
then $p=0$.  For $k=0$ it's obvious. We assume that 
 for any $p$ and $k$
$$\widetilde{V}^{(p)}\subset V^{(k)}+[V^{(k)},\mathfrak{B}]\; \; \; *$$

and  prove that:
 $$\widetilde{V}^{(p')}\subset V^{(k+1)}+[V^{(k+1)},\mathfrak{B}].$$
We have $\forall k, [V^{(k)},V^{(k)}]\subset [V^{(0)},V^{(k)}]$
and $ [V^{(k)},[V^{(k)},\mathfrak{B}]\subset (V^{(k)},V^{(k)},\mathfrak{B})$
since by definition $V^{(k+1)}= (V^{(k)},V^{(k)},\mathfrak{B}).$ So by
taking the commutator  of the relation (*) we have
$$ [\widetilde{V}^{(p)},V^{(k)}]\subset  [V^{(k)},V^{(k)}]+[[\mathfrak{B},V^{(k)}],V^{(k)}]$$
$$\;\; \; \; \; \; \; \; \; \; \; \; \; \; \; \; \subset  V^{(k+1)}+[[\mathfrak{B},V^{(k+1)}].$$

Hence the result.

Poof of the theorem.

Since the power of the right side of the theorem is a function of $k$, then by 
applying the lemma up to the reduction, we obtain the result. Hence
the theorem is proof.

From this theorem we have the following corollary:

{\bf Corollary} If $\mathfrak{B}$ is a solvable Bol algebra then 
the universal enveloping Lie algebra 
$\mathfrak{G}=\mathfrak{B} \dotplus \mathfrak{B}\wedge \mathfrak{B}$
is solvable.

\section{The Killing-Ricci form for Bol algebras and semisimple Bol algebras}

In this section we shall follow the construction in \cite{kik3}. A  Lie 
triple algebras (L.t.a) is  defined, as an anticommutative algebra
over a field $F$, whose multiplication is denoted by $XY$ for $X,Y \in \mathfrak{U}$.
Denote by $\mathbb{D}(X,Y)Z$ the inner derivative satisfying a series of 
identities. Bol algebras have the same construction as L.t.a; only that the 
endomorphism  $\mathbb{D}(X,Y)Z$ is called the inner pseudo-differential of the 
algebra $\mathfrak{B}$. The universal enveloping Lie algebra of $\mathfrak{B}$ 
is a Lie algebra $\mathfrak{G}=\mathfrak{B} \dotplus \mathbb{D}(\mathfrak{B},\mathfrak{B})$
such that the commutator in $\mathfrak{G}$ is defined as: 
$ [\xi, \eta]=\xi \cdot \eta + \mathbb(\xi,\eta)$. Conversely if we have any Lie
 algebra $\mathfrak{G}$ and a subalgebra $\mathfrak{h}$, we can define the 
operation in  $\mathfrak{B}$ such that  $\xi \cdot \eta=\prod_{\mathfrak{B}}[\xi,\eta]$.
The projection, is parallel to  $\mathfrak{h}$. As we have seen that Bol algebra allow embedding in Lie algebra, let  $\mathfrak{B}$ be a Bol algebra  let
$e_i, i=\overline{1,n}$ be its n-dimensional basis and 
$\mathbb{D}_1,......\mathbb{D}_N$
the basis of the  pseudo-derivation space $\mathbb{D}(\mathfrak{B},\mathfrak{B})$.
We assume that $\mathbb{D}(\mathfrak{B},\mathfrak{B}) \neq 0$ we define the 
operation in Bol algebra by means of a tensor as:
$$ e_i \cdot e_j=T_{i,j}^k e_k; \mathbb{D}(e_i,e_j)e_k=R_{i,j,k}^l e_l$$
$$ \mathbb{D}(e_i,e_j)=\mathbb{D}_{i,j}^{\tau}\mathbb{D}_{\tau}$$
$$[\mathbb{D}_{\tau},e_i]=\mathbb{D}_{\tau}e_i=K_{\tau,i}^j e_j$$
where the indices $i,j,k,l \in \{1,....n\}$ and $\tau \in \{1,......,N\}$
and 
$$\mathfrak{G}=\mathfrak{B} \dotplus \mathfrak{B}\wedge \mathfrak{B}$$
its universal enveloping Bol algebra. We denote by $\alpha$, the Killing
form of the universal enveloping Lie algebra. By the killing-Ricci form
$\beta$ of the Bol algebra we understand a symmetric bilinear form on
$\mathfrak{B}$ determined by the restriction of $\alpha$ to
 $\mathfrak{B}\wedge \mathfrak{B}$.

{\bf Proposition 1.} For any $\xi, \eta \in \mathfrak{B}$, we have the following
relation:
$$ \beta(\xi.\eta)=tr(L(\xi,\eta)+L(\eta,\xi))$$ where $L(\xi,\eta)\zeta=(\xi,\eta,\zeta).$

Proof 

Let $e_1,.....e_n$ be a basis on $\mathfrak{B}$ and 
$\mathbb{D}_1,.......,\mathbb{D}_N$ basis in the  pseudo-derivation  algebra.
We now calculate the basis element of the form $\beta$.

$$\beta(e_i,e_j)=tr(L_i,L_j)=tr(L_i L_j).$$
We also have $$ L_i L_j e_k =[e_i,[e_j,e_k]]=R_{i,j,k}^l e_l$$ and

$$L_i L_j \mathbb{D}_{\tau}=[e_i,[e_j,\mathbb{D}_{\tau}]]=$$
$$\; \; \; \; \; [e_i,K_{j\tau}^m e_m=K_{j\tau}^m T_{im}^{\nu} e_{\nu}.$$

Therefore $ \beta(e_i,e_j)=R_{ijk}^l +K_{j\tau}^m T_{im}^{\tau}$.
 On the other hand, the expression
$$[e_j,[e_i,e_m]]=T_{im}^{\rho}K_{j\rho}^s e_s=R_{jim}^s e_s $$
imply:
$$ R_{jim}^s=T_{im}^{\rho}K_{i\rho}^s, R_{ji}^m=T_{im}^{\rho}K_{j\rho}^m.$$
Therefore the basic element of the form $\beta$ can be written as:
$$ \beta(e_i,e_j)=R_{ijk}^k+R_{jik}^k.$$
Considering the Ricci tensor by convolution the upper and the lower index we 
then obtain:
$$R_{ij}=R_{ijk}^k=trL(e_i,e_j)$$
$$ \beta(e_i,e_j)=R_{ij}+R_{ji}=tr(L(e_i,e_j)+L(e_j,e_i))$$

Hence the proposition is proved.

An elementary check up shows that $\beta$ induces the properties of the form 
$\alpha$ for the Lie algebra $\mathfrak{G}$.

{\bf Definition 7.} We say that a Bol algebra $\mathfrak{B}$ is semisimple if it
contains only trivial ideals.

Let $\mathfrak{b}$ be an invariant symmetric bilinear form on $\mathfrak{B}$
satisfying the following condition:

$$\mathfrak{b}(X\cdot Y,Z)=\mathfrak{b}(X,Y\cdot Z)\; \; \; \; \; (**)$$
$$\mathfrak{b}((X, Y,Z),t)=\mathfrak{b}(Z,(X,Y,t))$$

Following \cite{kik3} we will introduce the notion of perpendicularity by means 
of the form $\mathbb{b}$ of an an object relatively to $\mathfrak{B}$.

Let $I$ be an ideal of the Bol algebra $\mathfrak{B}$, in particular we can 
assume $I \triangle \mathfrak{B}$. We will define the following two sets i.e 
the
 left orthogonal and the right orthogonal respectively by:

$$ I_l^{\perp}=\{X\in \mathfrak{B}/ \mathfrak{b}(X,I)=0\}$$
$$ I_r^{\perp}=\{X\in \mathfrak{B}/ \mathfrak{b}(I,X)=0\}$$

After introducing these concepts, this leads us to  the following 
definition:

{\bf Definition 8.} The set $\mathfrak{C}$ defined as
$$\mathfrak{C}=\{X\in \mathfrak{B}/ \mathfrak{B}\cdot X=0, (\mathfrak{B},\mathfrak{B},X)=0\}$$ it is called the center.

{\bf Proposition 2.} $\mathfrak{C}_r^{\top}=\mathfrak{C}_l^{\top}=\mathfrak{B}^{(1)}$
where $\mathfrak{B}^{(1)}=\mathfrak{B} \cdot \mathfrak{B}$.

Proof

Let $ X\in (\mathfrak{B}^{(1)})^{\top}_l $  This implies that $ \mathfrak{b}(X,\mathfrak{B},\mathfrak{B})=0$.
By using $(**)$ we have  $ \mathfrak{b}(X,(\mathfrak{B},\mathfrak{B},\mathfrak{B}))=0$ and hence  $ (X\cdot \mathfrak{B},\mathfrak{B})=0$. Also  
$ \mathfrak{b}((\mathfrak{B},\mathfrak{B},X),\mathfrak{B})=0$ gives  
$X\cdot \mathfrak{B}=0$ and  $(\mathfrak{B},\mathfrak{B},X)=0$ hence, $X \in C$
.
 Therefore,  $ X\in (\mathfrak{B}^{(1)})^{\top}_l \subseteq C$. In the same way, we
 prove the converse and therefore  $ X\in (\mathfrak{B}^{(1)})^{\top}_l=C $.
Analogously; one can prove that  $ X\in (\mathfrak{B}^{(1)})^{\top}_r=C $.
 Hence by taking the orthogonality twice we obtain the result

 $\mathfrak{C}_r^{\top}=\mathfrak{C}_l^{\top}=\mathfrak{B}^{(1)}$

Following Lie algebra theory and the result in \cite{kik3} we can state the 
following Theorem:

{\bf Theorem 3.} Let $\mathfrak{B}$ be a finite dimensional Bol algebra and 
let 
 $\mathfrak{b}$ be an invariant non-degenerated form on $\mathfrak{B}$ that 
does not contain  trivial ideal of the form:
$$ I^{(1)}=I\cdot I+(I,I,\mathfrak{B})=0$$
then, $\mathfrak{B}$ can be split in a direct sum of ideal
$$\mathfrak{B}=\mathfrak{B}_1+\mathfrak{B}_2+.......+\mathfrak{B}_r$$
where $b(\mathfrak{B}_i,\mathfrak{B}_j)=0$ for $i \neq j$ and each of the ideal
$\mathfrak{B}_i$ is a simple Bol algebra.

Proof

Let $I$ be an ideal of $\mathfrak{B}$. We can define the sets $I^{\top}_r$ and
$I \cap I^{\top}_r$  are also ideals of $\mathfrak{B}$. Then denoting by
$M=I \cap I^{\top}_r$ which is normal in $\mathfrak{B}$ and we will have
$b(M,M)=0,b(\mathfrak{B},M\cdot M)=0$;
$ b(\mathfrak{B},(M,M,\mathfrak{B}))=b(\mathfrak{B},(M,\mathfrak{B},M))=0 $.

Hence $M^{(1)}=0$ and $M=0$. We then have $\mathfrak{B}=I \dotplus I^{\top}_r $,
and $b(I, I^{\perp}_r)=0$. Now we must  prove  $b(I^{\perp}_r,I)=0$. We
have $I \subset \mathfrak{B}$ an ideal and so we can define
$ I^{(1)}=I\cdot I+(I,I,I+I^{\top}_r)\neq 0$ which is also an ideal of $I$ 
therefore $I=I^{(1)}$ and $I$ is a semisimple Bol algebra and  $b(I, I^{\perp}_r)=0$.
 Since $I\cdot I^{\top}_r$ and $(I,I,I^{\top}_r)\subseteq I \cap I^{\top}_r=0$
then $b(I^{\top}_r,I\cdot I)=b(I\cdot I^{\top}_r,I)=0$ and
 $b(I^{\top}_r,(I,I,I))=b((I,I,I^{\top}_r),I)=0$ also
 $b(I^{\top}_r,I^{(1)})=b( I^{\top}_r,I)=0$ therefore $I^{\top}_r \subset I^{\top}_l$
now the rest of the theorem is established by the induction on $dim\mathfrak{B}$.

As in \cite{kik3} we can state the following theorem for Bol algebra

{\bf Theorem 4.} Let $\mathfrak{B}$ be a finite dimensional Bol algebra over a field
of characteristic zero and let $\beta$ be the Killing-Ricci form then:
\begin{enumerate}
\item The universal enveloping Lie algebra $\mathfrak{G}$ of $\mathfrak{B}$ is
solvable if and only if $\mathfrak{B}$ and $(\mathfrak{B} ,\mathfrak{B} ,\mathfrak{B})$ are orthogonal.
\item  $\mathfrak{G}$ is simple if and only if $\beta$ is non-degenerated.
\item If the form $\beta$ is non degenerate and invariant, then  $\mathfrak{B}$
can be split into the direct sum of simple ideal which are orthogonal in the
sense of the definition above.

$$\mathfrak{B}=\mathfrak{B}_2\dotplus \mathfrak{B}_2+....... +\mathfrak{B}_r $$

and, the universal enveloping Lie algebra also split  into a direct sum of 
ideal
 $\mathfrak{G}_i$ such that each  $\mathfrak{G}_i$ is a universal enveloping
simple Lie algebra and we have 
 $\mathfrak{B}=( \mathfrak{B},\mathfrak{B},\mathfrak{B}).$
\end{enumerate}

Proof

In the proof of this theorem, we have to consider all the results of Lie 
algebras relatively to the Killing form.

\begin{enumerate}
\item  If  $\mathfrak{G}$ is a solvable Lie algebra then
$\alpha ( \mathfrak{G}, [\mathfrak{G},\mathfrak{G}])=0$.  As we know,
 $[\mathfrak{G},\mathfrak{G}]=(\mathfrak{B}, \mathfrak{B}, \mathfrak{B})\dotplus [\mathfrak{B}, \mathfrak{B}]$
and so 
 $\beta (\mathfrak{B}, (\mathfrak{B},\mathfrak{B},\mathfrak{B}))=0$.

conversely if $\beta (\mathfrak{B}, (\mathfrak{B},\mathfrak{B},\mathfrak{B}))=0$
then 
 $ [\mathfrak{B}, \mathfrak{B}]=Der\mathfrak{B} $ (set of pseudo-derivation)
and  so 
$\alpha(\mathfrak{B}, [\mathfrak{G}, \mathfrak{G}])=0$. It follows that
$\mathfrak{B}$ lies in the radical of $\mathfrak{G}$ and hence
$\mathfrak{B}$ is solvable. Moreover,  using the result of the corollary above 
about solvability we say the enveloping Lie algebra is also solvable.
\item Let $ \mathfrak{G}$ be a simple Lie algebra then the Killing form on
$\mathfrak{G}$ is non degenerate Since $\mathfrak{G}$ can be split into a 
direct sum of orthogonal subspace $\mathfrak{B}$ and
 $ [\mathfrak{B}, \mathfrak{B}]$
relatively to $\alpha$, then these subspaces are also non degenerate.
In particular, the form $\beta$ is non degenerate. Conversely, if the form
$\beta$ is non degenerate, then we have a natural embedding $\sigma$ see \cite{lis} on $\mathfrak{G}$:

$$\mathfrak{G}^{\sigma}=\mathfrak{B}^{\sigma}+[\mathfrak{B}^{\sigma},\mathfrak{B}^{\sigma}]$$ 

such that $X \longrightarrow X^{\sigma}$.

The operation in $\mathfrak{G}^{\sigma}$ using the properties of the
pseudo-derivative in $\mathfrak{B}$ is then  defined as:

$$[X^{\sigma},X^{\sigma}]=D(X,Y)$$
$$[X^{\sigma},D(Y,Z)]=(X,Y,Z)^{\sigma}$$
$$[D(X,Y),D(U,V)]=D(X,Y)D(U,V)-D(U,V)D(X,Y)$$
$$\; \; \; \; \; \;=D((X,U,V),Y)+D(X,(Y,U,V))$$

where $X,Y,Z,U,V \in \mathfrak{B}$.

Let $\alpha'$ be the Killing form in  $\mathfrak{G}^{\sigma}$, $\beta'$ the
restriction of  $\alpha'$ in  $\mathfrak{B}^{\sigma}$ then
$$\beta'(X^{\sigma},X^{\sigma})=\beta(X,Y)$$

since $\beta'$ is non degenerate on  $\mathfrak{B}^{\sigma}$, one can verify
easily that $\alpha'$ is also non degenerate; it is clear that the
subspace  $\mathfrak{B}^{\sigma}$ and $[\mathfrak{B}^{\sigma},\mathfrak{B}^{\sigma}]$
are orthogonal relatively to $\alpha'$. So by a simple verification one can 
see that  $[\mathfrak{B}^{\sigma},\mathfrak{B}^{\sigma}]$ is also non degenerate. Indeed if $\mathfrak{D} \in [\mathfrak{B}^{\sigma},\mathfrak{B}^{\sigma}]$
then  $\alpha'([\mathfrak{B}^{\sigma},\mathfrak{B}^{\sigma}],\mathfrak{D})=0$.
On the other hand,

 $$\alpha'([X^{\sigma},Y^{\sigma}],\mathfrak{D})=\beta'((X^{\sigma},(Y\mathfrak{D})^{\sigma})$$ 

$\forall X,Y \in \mathfrak{B}$ then $ (Y\mathfrak{D})^{\sigma}=0$. That means
$Y\mathfrak{D}=0$ and hence $\mathfrak{D}=0$. Therefore $\alpha'$ non degenerated 
and $\mathfrak{G}^{\sigma}$ semisimple hence the trilinear operation is
semisimple.
\item  Let us  assume that $\beta$ is non degenerate and invariant as in the
second point; we have the trilinear operation semisimple, and so
$\mathfrak{B}$ is semisimple. By using Theorem 2. which states that 
$\mathfrak{B}$ can be split into a direct sum of simple orthogonal ideal
$\mathfrak{B}_i$ relatively to $i=\overline{1,n}$.

\; \; Indeed if $X_i \in \mathfrak{B}_i$ and   $X_j \in \mathfrak{B}_j$ 
and $i \neq j$ then $X_i \cdot Y_j=0$ and 
$(t,X_i,Y_j)\in  \mathfrak{B}_i \cap \mathfrak{B}_j=0$ $\forall t \in \mathfrak{B}$
therefore $\mathfrak{D}(X_i,Y_j)$ is reduced to the null vector which is the
subset of $\mathfrak{B}$ With zero component $X_i Y_j$ hence we have

$$\mathfrak{G}=\mathfrak{G}_1+ \mathfrak{G}_2+\mathfrak{G}_3+....+\mathfrak{G}_n$$ 
where $\mathfrak{G}_i$ is the enveloping algebra of $\mathfrak{B}_i$
hence the theorem is proved.
\end{enumerate}

{\bf Aknoledgment:} This paper was able to be achieved, thanks to the help
of Professor M.  Kikkawa who sent  all the documentation I needed. And the 
good facility created at the ICTP Trieste, Italy.I would also like to express my thanks to professor Nguiffo Boyom for his valuable suggestions.

\end{document}